\documentclass[12pt,oneside]{amsart}
\usepackage{tkz-euclide}
\usepackage{geometry}                
\usepackage{graphicx}
\usepackage{amssymb}
\usepackage{mathtools}

\usepackage{amsmath,amsthm,amscd,amssymb}
\usepackage{latexsym}
\usepackage[colorlinks,citecolor=red,pagebackref,hypertexnames=false]{hyperref}
\usepackage{ulem}
\usepackage{soul}
\numberwithin{equation}{section}
\theoremstyle{plain}
\newtheorem{theorem}{Theorem}[section]

\newtheorem{corollary}[theorem]{Corollary}
\newtheorem{proposition}[theorem]{Proposition}

\theoremstyle{definition}
\newtheorem{definition}[theorem]{Definition}

\theoremstyle{remark}
\newtheorem{remark}[theorem]{Remark}

\newtheorem{case[theorem]}{Case}

\def\blue{\textcolor{blue}}

\def\N{\mathbb N}
\def\R{\mathbb R}
\def\hd{{\dim_{\mathcal H}}}
\def\vt{\vec{t}}
\def\supp{\hbox{supp}\,}

\date{\blue{April 30, 2024}}      
\author{P. Bhowmik, A. Greenleaf, A. Iosevich, S. Mkrtchyan, and F. Rakhmonov}
\address{Department of Mathematics, University of Rochester, Rochester, NY}
\email{pbhowmik@ur.rochester.edu}

\address{Department of Mathematics, University of Rochester, Rochester, NY}
\email{allan@math.rochester.edu}

\address{Department of Mathematics, University of Rochester, Rochester, NY}
\email{iosevich@math.rochester.edu}

\address{Department of Mathematics, University of Rochester, Rochester, NY}
\email{sevak.mkrtchyan@rochester.edu}

\address{Department of Mathematics, University of Rochester, Rochester, NY}
\email{frakhmon@ur.rochester.edu}

\thanks{AG is supported in part by National Science Foundation award DMS-2204943.
AI was supported in part by the National Science Foundation under grant no. HDR TRIPODS - 1934962 
and by the NSF DMS - 2154232. 
SM is supported in part by the Simons Foundation Grant no. 523555.}

\begin{document}

\title{Similar  point configurations via group actions}

\maketitle

\begin{abstract} We prove that for $d\ge 2,\, k\ge 2$, if the Hausdorff dimension of a compact set 
$E\subset \mathbb{R}^d$ is greater than $\frac{d^2}{2d-1}$, then, for any given $r > 0$, 
there exist $(x^1, \dots, x^{k+1})\in E^{k+1}$, $(y^1, \dots, y^{k+1})\in E^{k+1}$, 
a rotation $\theta \in \mathrm{O}_d(\mathbb{R})$, and a vector $a \in \mathbb{R}^d$ 
such that $rx^j = \theta y^j - a$ for $1 \leq j \leq k+1$. 
Such a result on existence of similar $k$-simplices in thin sets had previously been established under a more stringent 
dimensional threshold in Greenleaf, Iosevich and Mkrtchyan \cite{GIM21}. The argument we are use to prove the main result here was previously employed
in Bhowmik and Rakhmonov \cite{BR23} to establish a finite field version. We also show the existence of multi-similarities of arbitrary multiplicity in $\R^d$, show how to extend these results from similarities to arbitrary proper continuous maps, as well as
explore a general  group-theoretic formulation of this problem in vector spaces over finite fields.
\end{abstract}  

\maketitle

\section{Introduction} 

\vskip.125in 

The celebrated Falconer distance conjecture (see e.g. \cite{Fal86}, \cite{M95}) states that if the Hausdorff dimension 
$\hd(E)$ of a compact set $E \subset \mathbb{R}^d$, $d \ge 2$, is greater than $\frac{d}{2}$, then the Lebesgue measure of the 
distance set $\Delta(E):=\{|x-y|: x,y \in E \}$ is positive. Currently,  the best  threshold known for the Falconer distance problem in 
two dimensions is $\frac{5}{4}$, due to Guth, Iosevich, Ou and Wang (\cite{GIOW20}). In higher dimensions, the best exponent 
known is $\frac{d}{2}+\frac{1}{4}-\frac{1}{8d+4}$, due to Du, Ou, Ren, and Zhang (\cite{DORZ23}). 

\smallskip

More generally, given a compact set $E\subset \mathbb{R}^d$, $d\geq 2$, let 
\begin{equation*}
   \Delta_k(E)=\{v(x^1, \dots, x^{k+1}): x^j \in E\}\subset \R^{k+1\choose 2}, 
\end{equation*}
where $v(x^1, \dots, x^{k+1})$ is the $\binom{k+1}{2}$-vector whose entries are $|x^i-x^j|$, $1 \leq i<j \leq k+1$, 
listed in the dictionary order. 

\smallskip

We can think of $\Delta_k(E)$   as the set of non-congruent $(k+1)$-tuples of elements of $E$
 (modulo the action of the symmetric group $S_{k+1}$). 
 More precisely, we say that $(x^1, \dots, x^{k+1})\in E^{k+1}$ is {\it congruent} to $(y^1, \dots, y^{k+1})\in E^{k+1}$ if there exists an 
 orthogonal matrix $\theta\in \mathrm{O}_d(\mathbb{R})$ and $a \in {\mathbb R}^d$ 
 such that $y^j=\theta x^j+a$ for $1 \leq j \leq k+1$. 
Equivalently, for $k\le d$, one can consider $\Delta_k(E)$ as the set
of congruence classes of $k$-simplices generated by points in $E$. 
\smallskip

One can also consider {\it similarity} classes of $k$-simplices.
Greenleaf, Iosevich, and Mkrtchyan  proved in  \cite{GIM21} that if the Hausdorff dimension of $E$, $\hd(E)$,  satisfies 
$\hd(E)>s_k:=\frac{dk+1}{k+1}$,  then

\begin{eqnarray} \label{weakconclusion} 
& & \hbox{ for every $r>0$, there exist distinct} \nonumber \\
 & & (x^1, \dots, x^{k+1})\in E^{k+1},\quad  
(y^1, \dots, y^{k+1})\in E^{k+1}\\
& & \hbox{  such that } |y^i-y^j|=r|x^i-x^j|  \hbox{ for every } 1 \leq i<j \leq k+1. \nonumber
\end{eqnarray} 

In other words, given any $r>0$, 
there exist two $(k+1)$-tuples of elements of $E$ which are similar via a translation, rotation, and the prescribed scaling factor $r$. 
More precisely, the following was proved in  \cite{GIM21}: 
Given $E$ with $\hd(E)>s_k$, let $s_k<s<\hd(E)$ and $\mu$ be a Frostman measure on $\R^d$, supported on $E$ 
and of finite $s$-energy, as in Frostman's Lemma (see \cite{M95}).
From $\mu$ one forms  the natural configuration measure   $\nu_k$  on $\Delta_k(E)$, 
i.e., for any $f\in C_0\left(\R^{k+1 \choose 2}\right)$,
$$ \int f(\vec{t} ) \, d\nu_k(\vec{t})=\int \dots \int f\left(v(x^1, \dots, x^{k+1})\right) d\mu(x^1)d\mu(x^2) \dots d\mu(x^{k+1}).$$ 
\smallskip

For $0<r<\infty$, define 
\begin{equation}\label{eqn Delta r}
    \Delta^r_k(E)=\left\{\vec{t} \in \Delta_k(E): r \vec{t} \in \Delta_k(E)\right\}\subset \Delta_k(E),  
\end{equation}
which records the $(k+1)$-point configurations in $E$ which also have an $r$-scaled similar copy in $E$.
Then, in  \cite{GIM21} it was shown that if $\hd(E)>s_k$, then
\begin{equation} \label{GIM21mainestimate}
 \nu_k(\Delta_k^r(E))>0, 
\end{equation}
which implies (and is in fact much stronger than) statement \eqref{weakconclusion}.

\medskip

In the first results of the current paper, Theorem \ref{main} and Corollary \ref{cor main},
we will show that while (\ref{GIM21mainestimate}) might require the full strength of the 
machinery used in \cite{GIM21}, one can derive the conclusion (\ref{weakconclusion}) for any $k$ under a much lower, $k$-independent 
dimensional threshold,
$\hd(E)>\frac{d^2}{2d-1}$, and extend it to $k>d$ (for which $(k+1)$-point sets no longer correspond to $k$-simplices). 
In brief, the work in \cite{GIM21} shows, and uses, that
if  $\hd(E)>s_k=\frac{dk+1}{k+1}$ defined above, then
the measure $\nu_k$ on 
$\R^{k+1 \choose 2}$, supported on $\Delta_k(E)$,
 is absolutely continuous with respect to Lebesgue measure $d\vec{t}$, 
with a density  in $L^2\left(\R^{k+1 \choose 2}\right)$.
In contrast  the method we employ in this paper 
merely requires  use of the  lower threshold, $\hd(E)>\frac{d^2}{2d-1}$,
 above which  the configuration measure $\nu_1$ on the distance set $\Delta(E)=\Delta_1(E)$ is
 known to be in $L^2({\mathbb R})$.
Theorem \ref{thm similar}  then uses a continuous version of the pigeonhole principle (Proposition \ref{prop pigeon}, from \cite{GIM21}) to 
extend these results to multi-similarities of any multiplicity, while Theorem \ref{thm transform} in Section \ref{transformations} shows how these results can be extended to more general transformations in $\R^d$.
\medskip

In addition to the Euclidean setting described so far, we also examine an analogous problem 
in the context of general transitive  group actions on vector spaces over 
finite fields. 
See Theorems \ref{BR23} and \ref{maingroupactiontheorem}, and Corollary \ref{sl2actioncorollary}
and their proofs in Section \ref{sec proof of 1.4}.
Our goal is to lay the groundwork for future explorations in this broader context, 
including   both finite field and Euclidean problems. 

\medskip

\subsection{Results in  Euclidean spaces} 
In the Euclidean setting, our main result is: 
\medskip

\begin{theorem} \label{main}
Let $E \subset \mathbb{R}^d$, $d \ge 2$, be compact with Hausdorff dimension satisfying $\hd(E)> \frac{d^2}{2d-1}$.
Then, for any $r>0$, there exist $\theta \in \mathrm{O}_d(\mathbb{R})$ and $a \in \mathbb{R}^d$ 
such that  $rE \cap (\theta E-a)$ has positive Hausdorff dimension, and is thus uncountably infinite.
In particular, for every $k\in\N$, \eqref{weakconclusion} holds.
\end{theorem}

\vskip.125in

An immediate consequence of Theorem \ref{main} is the following.
\smallskip

\begin{corollary}\label{cor main}
Let $E \subset \mathbb{R}^d$, $d \ge 2$, be a compact set with
$\hd(E)>\frac{d^2}{2d-1}$.
Then, for any $k \ge 1$ and $r > 0$, 
there exist distinct $(x^1, \dots, x^{k+1})\in E^{k+1}$ and $(y^1, \dots, y^{k+1})\in E^{k+1}$, 
such that $|x^i-x^j| = r|y^i-y^j|$ for all $1 \leq i < j \leq k+1$. 
Thus, if $1\le k\le d$, for any $r>0$ there exist $r$-similar $k$-simplices in $E$.
\end{corollary}
\smallskip

As in \cite{GIM21} for the higher threshold $s_k$, a measure-theoretic pigeonhole principle allows one to leverage 
the proof of Theorem \ref{main} to 
obtain the existence of {\it multi-similarities} in $E$, of arbitrarily large multiplicity, for $\hd(E)$ above the new, lower threshold:

\begin{theorem} \label{thm similar} 
Let $E \subset \mathbb{R}^d$, $d \ge 2$, be a compact set with
$\hd(E)>\frac{d^2}{2d-1}$.
Then for all $n\ge 1$
and any compact interval $I\subset (0,\infty)$, 
there exists  $M=M(n,E,I)\in\N$ such that for any  distinct $r_1,\dots,r_{M}\in I$, there
exist distinct $r_{i_1},\dots, r_{i_n}$ such that: For every $k\ge 1$, there is a
$\vec{t}\in\R^{{k+1}\choose 2}\, $such that
$$\left\{r_{i_1}\vt,\, r_{i_2}\vt,\, \dots\, ,\, r_{i_n}\vt\, \right\}\subset \Delta_k(E),$$
i.e., there exists an $n$-similarity of $(k+1)$-point sets in $E$ with scaling factors coming from among $\{r_i\}_{i=1}^M$.
\end{theorem}

See Section \ref{subsec multi} for the relevant definitions, statements and proof.
\medskip

We also extend Theorem \ref{main} by showing that the dilations in Theorem \ref{main} 
can be replaced by more general  transformations:

\begin{theorem}\label{thm transform}
Suppose $T:\mathbb{R}^d\to\mathbb{R}^d$ is a proper continuous  map,
and $E\subset\R^d$ is compact. Given $s<\hd(E)$, let $\mu$ 
be a Frostman measure $\mu$ supported on $E$ and of finite $s$-energy. 
Define the pushforward measure $\mu_T$ by 
    its action on $f\in C_0(\R^d)$,
    $$ \int f(x) d\mu_T(x):=\int f(T(x)) d\mu(x).$$
    If $s$ is such that
    $$\int {|\widehat{\mu}_T(\xi)|}^2 {|\xi|}^{-s \frac{d-1}{d}} d\xi \, <\, \infty
    $$
 then $T(E) \cap (\theta E-a)$ has positive Hausdorff dimension, and is thus uncountably infinite.
 Thus, for any $k\in\N$,
 we can find points $(x_1,\dots,x_{k+1})\in E^{k+1}$ and $(y_1,\dots,y_{k+1})\in E^{k+1}$  such that $T(x_i)=\theta y_i-a$ for all $i$ for some $\theta \in \mathrm{O}_d(\mathbb{R})$.
 \end{theorem}
 
 For example, Theorem \ref{thm transform} applies  if $T$  is a diffeomorphism of $\R^d$ 
 which is the identity outside of a compact set, since such maps preserve Sobolev spaces of all orders.

\medskip

\subsection{Finite field setting and transitive group actions}

The idea behind the proof of \eqref{weakconclusion} was originally implemented in the finite field setting by 
Bhowmik and Rakhmonov  in \cite{BR23}. 
\smallskip

Let $\mathbb{F}_q^d$ be the  vector space of dimension $d\ge 2$  over the finite field $\mathbb{F}_q$ with $q$ elements. 
We define a function $\lVert \cdot \rVert: \mathbb{F}_q^d \to \mathbb{F}_q$ 
by $\lVert \alpha \rVert \coloneqq \alpha_1^2 + \dots + \alpha_d^2$ for $\alpha = (\alpha_1, \dots, \alpha_d) \in \mathbb{F}_q^d$. 
Note that this is not a norm, as it is $\mathbb F_q$-valued (and homogeneous of degree 2), and we do not impose any metric 
structure on $\mathbb{F}_q^d$. 
However, $\lVert \cdot \rVert$ does share an important feature of the Euclidean norm: it is invariant under orthogonal 
transformations.
\smallskip

For $m\ge 2$, let $(\mathbb{F}_q)^m\coloneqq\{a^m: a\in \mathbb{F}_q\}$;
in particular, $(\mathbb{F}_q)^2$  is the set of quadratic residues in $\mathbb{F}_q$. 
Then the following holds.

\smallskip
\begin{theorem} \label{BR23} \cite[~Thm. 1.3]{BR23}
Let $k\ge 1$.
Suppose $\emptyset \neq A\subset \{(i,j):1\leq i<j\leq k+1\}$
and $r\in (\mathbb{F}_q)^2\setminus \{0\}$. 
If $E\subset \mathbb{F}_q^d$ with $|E|\geq \sqrt{k+1}q^{\frac{d}{2}}$, 
then there exist $(x_1,\ldots,x_{k+1})\in E^{k+1}$ and $(y_1,\ldots,y_{k+1})\in E^{k+1}$ 
such that $\lVert y_i-y_j \rVert=r\lVert x_i-x_j \rVert$ for $(i,j)\in A$, and $x_i\neq x_j$, $y_i\neq y_j$, for all $1\leq i<j\leq k+1$.
\end{theorem}
\smallskip

As a straightforward corollary, through the variation of the underlying set $A$, we can establish thresholds for the existence of dilated $k$-cycles, $k$-paths, and $k$-stars (for $k\geq 3$) with a dilation ratio $r\in (\mathbb{F}_q)^2\setminus \{0\}$. Notably, Rakhmonov has previously examined the cases of $2$-paths, $4$-cycles, and $k$-simplices in \cite{MR4609035}.

\smallskip

The  original idea for a proof of Theorem \ref{BR23} was based on the a group-theoretic approach. 
(However, in \cite{BR23}, the first and fifth listed authors opted for an alternative, weaker version of this method.) 
This group-theoretic approach can be generalized to arbitrary transitive group actions, 
which we state here and prove in  Section \ref{sec proof of 1.4}:

\smallskip
\begin{theorem} \label{maingroupactiontheorem} 
Suppose that $G$ is a finite group acting transitively on a set $X$, and let 
$E$ and $H$  be subsets of $X$. Then 
\begin{equation}\label{lower bound for the intersection}
    \max_{g\in G}|H\cap gE|\geq \frac{|H||E|}{|X|}.
\end{equation}
\end{theorem} 
\smallskip

Another application of  Theorem 
\ref{maingroupactiontheorem} pertains to $\textup{SL}_d({\mathbb F}_q)$ actions:

\smallskip
\begin{corollary} \label{sl2actioncorollary}
Suppose $k \ge d$.
Let $E \subset \mathbb{F}_q^d$, with $|E| \ge \sqrt{k+1}q^{\frac{d}{2}}$. 
Then for every $r\in \left(\mathbb F_q\right)^d\setminus \{0\}$,
there exist $(x_1, \dots, x_{k+1}) \in E^{k+1}$ and $(y_1, \dots, y_{k+1}) \in E^{k+1}$ such that 
$$\det(x_{i_1}, x_{i_2}, \dots, x_{i_d}) = r \det(y_{i_1}, y_{i_2}, \dots, y_{i_d})$$
for all $d$-tuples of elements from $(x_1, \dots, x_{k+1}) \in E^{k+1}$ and $(y_1, \dots, y_{k+1}) \in E^{k+1}$, respectively, where $\det(u_1, u_2, \dots, u_d)$ is the determinant of the $d \times d$ matrix where the columns are vectors in $\mathbb{F}_q^d$. 
\end{corollary}
\bigskip

The paper is organized as follows. In Section \ref{sec proof of main} we prove Theorem \ref{main} in $\R^d$,
and recall in Section \ref{subsec multi}
the  material needed to state  and prove 
Theorem \ref{thm similar} on multi-similarities.
Then, in Section \ref{transformations} we prove Theorem \ref{thm transform}, extending the context of Theorem \ref{main} 
from dilations
to more general transformations in $\R^d$.
Finally, in Section \ref{sec proof of 1.4} we prove the group action result  Theorem \ref{maingroupactiontheorem}
and  show how  it implies Theorem \ref{BR23}
and Corollary \ref{sl2actioncorollary} in the finite field setting.

\vskip.125in 

\section{Similarities and multi-similarities in \texorpdfstring{$\mathbb{R}^d$}{Rd}}\label{sec proof of main}

\subsection{Proof of Theorem \ref{main}}\label{subsec proof of main} 

\vskip.125in 

Let $E\subset\R^d$ be compact.
Then, for any $s<\hd(E)$, we can equip $E$ with a Frostman measure $\mu$, a probability measure supported on $E$, 
of finite $s$-energy and 
satisfying a dimension $s$  ball condition. 
For $0<r<\infty$ and $\theta\in O_d(\R)$, the group of orthogonal transformations of $\R^d$, 
define $\mu_r$ and $ \mu^\theta$, also probability measures of dimension $s$,   
by their actions on $f\in C_0(\R^d)$, 
$$ \int f(x) d\mu_r(x):=\int f(rx) d\mu(x),\quad \int f(x) d\mu^\theta(x):=\int f(\theta x) d\mu(x),$$
with Fourier transforms $\widehat{\mu_r}(\xi)=\widehat{\mu}(r\xi)$ and  
$\widehat{\mu^\theta}(\xi)=\widehat{\mu}\left(\theta^{-1}\xi\right)$.
Using this notation,    consider
\begin{equation}\label{eqn to justify}
 \int \left(\mu_r*\mu^\theta\right)(x) dx.
 \end{equation}
 If the expression in \eqref{eqn to justify} converges we are done because, formally,  it equals 
$$\left( \int d\mu_r(x)\right) \cdot \left( \int d\mu^\theta(x)\right) = \int d\mu_r(x) \cdot \int d\mu^\theta(x)=1.$$
However, in order to decouple the integrals, one needs to demonstrate that \eqref{eqn to justify} converges, 
at least for a  set  of $\theta$ of full measure. If this convergence is established, and if we define
 $$A_{r,\theta}:=\supp\left(\mu_r*\mu^\theta\right)\subset\R^d,$$
then one has
$$1= \int \left(\mu_r*\mu^\theta\right)(x) dx =\int{\mathbf 1}_{A_{r,\theta}}(x)\cdot \left(\mu_r*\mu^\theta\right)(x)\, dx.$$
Squaring this, applying Cauchy-Schwarz and then 
integrating  with respect to the normalized Haar measure $d\theta$ on $O_d(\R)$, 
 it follows that
\begin{eqnarray*}
1=1^2 &=& \int_{O_d(\R)}\left( \int_{\R^d}  {\mathbf 1}_{A_{r,\theta}}(x)\cdot \left(\mu_r*\mu^\theta\right)(x) dx\, \right)^2   d\theta \\
& \le & \int_{O_d(\R)} \left(\int_{\R^d} {\mathbf 1}_{A_{r,\theta}}(x)^2 dx \right)
\cdot \left( \int_{\R^d} \left(\mu_r*\mu^\theta\right)\!(x)^2 dx \right)\, d\theta \\
&=& \int {\mathcal L}^d \left(A_{r,\theta}\right) \cdot \int {(\mu_r*\mu^\theta)(x))}^2 dx \, d\theta   \\
&\leq& \left(\sup_{\theta\in O_d(\R)} {\mathcal L}^d\left(A_{r,\theta}\right)\right) \cdot \int \int {(\mu_r*\mu^\theta)(x))}^2 
dx\,  d\theta,    
\end{eqnarray*}
where $\mathcal L^d$ denotes Lebesgue measure on $\R^d$.
Hence,
\begin{equation}  \label{keysetup} 
\sup_{\theta} \mathcal{L}^d\left(A_{r,\theta}\right) \ge \frac{1}{\int \int (\mu_r*\mu(\theta \cdot)(x))^2  dx\,  d\theta}. 
\end{equation}

\medskip

To obtain a lower bound for the left hand side of \eqref{keysetup},
we need an upper bound on  the denominator $\int \int (\mu_r*\mu(\theta \cdot)(x))^2  dx\,  d\theta$ of the right hand side.
To obtain this, we  use the following result.

\begin{theorem} \label{DZtheorem}  (Wolff \cite{W99} for $d=2$; Du and Zhang \cite{DZ19} for $d\ge3$)
Let $\mu$ be a  measure, supported on a compact set  on $\mathbb{R}^d$
and of dimension $\alpha$  for some $\alpha\in\left(\frac{d}2,d\right)$.
Then,
$$ \int_{S^{d-1}} |\widehat{\mu}(R\omega)|^2 \, d\omega \leq C_{\alpha,{\mu}}R^{-\alpha \frac{d-1}{d}}. $$
\end{theorem}

\medskip

Returning to the proof  of Theorem \ref{main}, define $F_{r,\theta}(x) = \left(\mu_r * \mu^\theta\right) (x)$.
Then, 

$$
    \int \int F^2_{r,\theta}(x) d\theta \, dx  
 = \int \int {|\widehat{\mu}_r(\xi)|}^2 {|\widehat{\mu}(\theta \xi)|}^2 d\theta \, d\xi    
\lesssim \int {|\widehat{\mu}_r(\xi)|}^2 {|\xi|}^{-\alpha \frac{d-1}{d}} d\xi,    
$$
where for the inequality we have used Theorem \ref{DZtheorem} for any $\frac{d}2<\alpha<\hd(E)$. 
Taking $\alpha=s$, the dimension of $\mu$, then
 if $s > \frac{d^2}{2d-1}$,
 the denominator of the right side of \eqref{keysetup}  is bounded above if $\hd(E)> \frac{d^2}{2d-1}$,
and thus the left hand side is bounded below, say by $c_r>0$.  This also shows that in fact the 
integral in \eqref{eqn to justify} converges, thereby justifying the calculations above.
\medskip

(Note for use  in the next subsection that $c_r$ depends on $s$ and the
$s$-energy of $\mu_r$;
thus, as $r$ ranges over any compact interval $I\subset (0,\infty)$, $c_r$  is uniformly bounded away from 0.)
\medskip

We have shown that, for any $r > 0$, there exists $\theta \in \mathrm{O}_d(\mathbb{R})$ such that 
$$ \mathcal{L}^d\left(A_{r,\theta}\right) \ge c_r > 0, $$ or equivalently 

$$ \mathcal{L}^d\left(\{a\in\R^d\, :\, \mu\{y\in E\, :\, (\exists\, x\in E)\, (\, rx=\theta y-a\,)\, \}>0
\}\right) \ge c_r > 0.$$

Fixing any one such $a$, it follows that $\{y\in E\, :\, (\exists\, x\in E)\, rx=\theta y-a\}$ has positive $\mu$ measure, 
thus has positive Hausdorff dimension, and hence is an (uncountably) infinite set.
In particular,  for any $k \ge 1$, 
there exist $(x^1, \dots, x^{k+1}) \in E^{k+1}$, $(y^1, \dots, y^{k+1}) \in E^{k+1}$, $\theta \in \mathrm{O}_d(\mathbb{R})$, 
and $a \in \mathbb{R}^d$ such that $rx^j = \theta y^j - a$ for $1 \leq j \leq k+1$. 
This completes the proof of Theorem \ref{main}. 
\vskip.25in 

\subsection{Extension to multi-similarities: Proof of Theorem \ref{thm similar}}\label{subsec multi}

By the parenthetical comment above, 
\begin{equation}\label{eqn uniform}
\left(\forall\, I\subset (0,\infty), I\text{ compact}\right)\, (\exists\, c_I>0)\,(\forall\, r\in I)\, c_r\ge c_I \, .
\end{equation}
Modifying slightly the proof of Theorem \ref{main} and using a continuous version of the pigeonhole principle from \cite{GIM21}
will allow us to prove the existence of arbitrarily many simultaneous similarities in $E$.
We start with some background.
\medskip

The set $\Delta_k^r(E)$ defined in \eqref{eqn Delta r} records the pairs $\{\vt,r\vt\}\subset \Delta_k(E)$
of $(k+1)$-point sets in $E$  which are similar by an $r$-scaling and translation. 
Similarly, one can look at $n$-tuples of similar $(k+1)$-point sets.
Recall the following definition from \cite{GIM21}; since all line segments are similar to each other, this is only interesting for $k\ge 2$.

\begin{definition}\label{def multi}
Let $n\ge 1$. A  {\it multi-similarity of multiplicity $n$} of $(k+1)$-point sets in $E$ is a set of the form
$\{r_1\vt,\, r_2\vt,\, \dots ,\, r_{n}\vt\}\subset \Delta_k(E)$, for some $\vt\in\R^{{k+1}\choose 2}$ 
and  $r_1,\, r_2,\, \dots, r_{n}\in (0,\infty)$, pairwise distinct.
(If $k\le d$,  this is a  simultaneous similarity of $k$-simplices in $E$ under $n$ different scalings.)
\end{definition}

In \cite{GIM21}, it was shown that the version of our Theorem \ref{main} with the higher threshold $\hd(E)>s_k$ had
as a corollary a multi-similarity version, \cite[Thm. 1.7]{GIM21}. 
This was done by combining a multi-similarity version of \eqref{GIM21mainestimate}
with a measure-theoretic version  of the pigeonhole principle \cite[Lem. 2.1]{GIM21}. 
We can combine this reasoning with a modification of the proof of Theorem \ref{main} 
to obtain Theorem \ref{thm similar}, i.e., 
the existence, above the same  lower threshold of Theorem \ref{main}, of multi-similarities of all multiplicities.
\medskip

We start from the same expression \eqref{eqn to justify} as above, but now consider $\theta$ as a variable rather than a parameter.
Define $D_r=\supp_{x,\theta}\left(\mu_r*\mu^\theta\right)\subset\R^d\times O_d(\R)$.
Then, from
$$1=\int\int_{\R^d\times O_d(\R)} \left(\mu_r*\mu^\theta\right)(x)\, dx\, d\theta = 
\int\int_{\R^d\times O_d(\R)} {\mathbf 1}_{D_r}(x,\theta)\cdot \left(\mu_r*\mu^\theta\right)(x)\, dx\, d\theta, 
$$
a variation of the calculation in the previous section yields
$$1\le \left({\mathcal L}^d\times d\theta\right)(D_r)\cdot \int\int F_{r,\theta}(x,\theta)^2\, dx\, d\theta $$
and thus
\begin{equation}\label{eqn similar lower bound}
\left({\mathcal L}^d\times d\theta\right)(D_r)\ge c_r.
\end{equation}

Now recall the

 \begin{proposition}\label{prop pigeon}
 {\bf Measure-theoretic Pigeonhole Principle (\cite{GIM21})} \\
 Let $\mathcal X=(X,\mathcal M,\sigma)$ be a  measure space with $\sigma(X)<\infty$.
 For  $0<c<\sigma(X)$, define $\mathcal M_c=\{A\in\mathcal M: \sigma(A)\ge c\}$. 
 Then, for every $n\ge 1$, there exists an\linebreak $N=N(\mathcal X,c,n)$ such that: for 
 any collection $\{A_1,\dots,A_N\}\subset \mathcal M_c$ of cardinality $N$, 
 there is a subcollection $\{A_{i_1},\dots, A_{i_n}\}$ of cardinality $n$ such that $\sigma(A_{i_1}\cap\cdots\cap A_{i_n})>0$ 
 and hence 
 $A_{i_1}\cap\cdots\cap A_{i_n}\ne\emptyset$.
 \end{proposition}

From \eqref{eqn uniform}, \eqref{eqn similar lower bound}, we see that
$\{D_r\}_{r\in I}$ is a family of measurable subsets of the finite measure space 
$\overline{B}(0,R)\times O_d(\R)$, equipped with $\mathcal L^d \times d\theta$,
where 
$$R=2(\sup \{r\,:\, r\in I\}) \cdot (\sup\{|x|\, :\, x\in E\}),$$
with measures bounded below by $c_I$. Applying Proposition \ref{prop pigeon}, there exists an $M\in\N$ such that,
for any distinct $\{r_i\}_{i=1}^M\subset I$, there is a subset $\left\{r_{i_j}\right\}_{j=1}^n$ such that
$$\left(\mathcal L^d \times d\theta\right)\left(\bigcap_{j=1}^n D_{r_j}\right)>0.$$
Thus, there exists a $\theta_0\in O_d(\R)$ such that
$$\mathcal L^d \left(\bigcap_{j=1}^n D_{r_j} \cap \{\theta=\theta_0\}\right)>0.$$
Since $\cap_{j=1}^n D_{r_j} \cap \{\theta=\theta_0\}$ has positive Lebesgue measure, it is infinite and hence for any $k\in\N$ 
we can find $\{x^l\}_{l=1}^{k+1}$ and $\{y^l\}_{l=1}^{k+1}$ such that
$$r_{i_j}x^l=\theta_0y^l-a,\quad\forall\,  1\le j\le n,\, 1\le l\le k+1.$$
This finishes the proof of Theorem \ref{thm similar}.
Note that what was established is considerably stronger than the statement of the theorem.

\medskip
\section{Proof of Theorem \ref{thm transform}}
\label{transformations}

Since $T:\mathbb{R}^d\to\mathbb{R}^d$ is proper and  continuous,  
 the Frostman measure $\mu$ on $E$ defines  a measure $\mu_T$ by  
$$ \int f(x) d\mu_T(x):=\int f(T(x)) d\mu(x),\quad\forall\, f\in C_0(\R^d).$$
Suppose that
$$\int {|\widehat{\mu}_T(\xi)|}^2 {|\xi|}^{-s \frac{d-1}{d}}\, d\xi\,  <\, \infty.
$$
Arguing as in the proof of Theorem \ref{main},  if we set 
$$A_{T,\theta}:=\supp\left(\mu_T*\mu^\theta\right)\subset\R^d,$$
one sees that there exists a $c_T>0$ such that 
\begin{equation*} 
\sup_{\theta\in O_d(\R)} \mathcal{L}^d\left(A_{T,\theta}\right) \ge c_T,
\end{equation*}
from which it follows that there exists a $\theta_0\in O_d(\R)$ such that 
$$ \mathcal{L}^d\left(\{a\in\R^d\, :\, \mu\{y\in E\, :\, (\exists\, x\in E)\, (\, T(x)=\theta_0 y-a\,)\, \}>0
\}\right) \ge c_T > 0.$$
It   follows that for all $k\in\N$, there are points $a\in\R^d$ and  $x_1,\dots,x_{k+1},y_1,\dots,y_{k+1}\in E$ and such that $T(x_i)=\theta_0 y_i-a$ for all $1\le i\le k+1$.
\medskip

\begin{remark}
Similarly to how one can go from a single dilation factor in Theorem \ref{main}
to the  multi-similarities in Theorem \ref{thm similar},
one can boost from the one transformation $T$ in Theorem \ref{thm transform} to multiple transformations. 
Suppose that $\mathcal T$ is a family of proper, continuous maps  and $s$ is such that there is a uniform bound,
$$\int {|\widehat{\mu}_{T}(\xi)|}^2 {|\xi|}^{-s \frac{d-1}{d}} d\xi\, \le C,\quad \forall\, 1\le j\le m,\, T\in\mathcal T.
$$ 

Following  the proof of Theorem \ref{thm similar} in Section \ref{subsec multi},
one can show that for all $\, n\in\N$, there exists an $M\in\N$ such that
for any distinct $\{T_j\}_{j=1}^M\subset \mathcal T$, there exists a subset $\left\{T_{j_l}\right\}_{l=1}^n$,
 a $\theta_0\in O_d(\R),\, a\in\R^d$, and for all $k\in\N$, points $x_1,\dots,x_{k+1},y_1,\dots,y_{k+1}\in E$  such that $T_{j_l}(x_i)=\theta_0 y_i-a$ 
 for all $1\le i\le k+1,\, 1\le l\le n$.
\end{remark}

\bigskip

\section{Proofs of  Theorem \ref{maingroupactiontheorem}, Theorem \ref{BR23}
and Corollary \ref{sl2actioncorollary} } \label{sec proof of 1.4}

\vskip.125in 

We start with the proof of Theorem \ref{maingroupactiontheorem}, 
and then show how to apply it to obtain the similarity set
results Theorem \ref{BR23} and Corollary \ref{sl2actioncorollary}.
\medskip

\noindent{\bf Proof of  Theorem \ref{maingroupactiontheorem}.} Suppose that $G$ is a finite group acting transitively on a set $X$, 
and let $H$ and $E$ be subsets of $X$. For each $x \in X$, let $\mathcal{O}_x \coloneqq \{gx : g \in G\}$ 
be the orbit of $x$, and $\mathcal{S}_x \coloneqq \{g \in G : gx = x\}$ be the stabilizer of $x$. 
The orbit-stabilizer theorem tells us that $|G| = |\mathcal{O}_x| \cdot |\mathcal{S}_x|$. 
Since the action of $G$ on $X$ is transitive, then $\mathcal{O}_x = X$, and hence $|\mathcal{S}_x| = {|G|}/{|X|}$.

\smallskip

Let $\mathcal{S}_{xy} \coloneqq \{g \in G : gx = y\}$. It is not difficult to prove that $|\mathcal{S}_{xy}| = |\mathcal{S}_x|$. Indeed, 
consider the map $\phi: \mathcal{S}_{xy} \times \mathcal{S}_{x} \to \mathcal{S}_{xy}$ defined by $(h, g) \xmapsto{\phi} hg$. We 
observe that every $u \in \mathcal{S}_{xy}$ has exactly $|\mathcal{S}_{xy}|$ preimages, 
i.e., $|\phi^{-1}(\{u\})| = |\mathcal{S}_{xy}|$. This immediately implies that 
$|\mathcal{S}_{xy} \times \mathcal{S}_{x}| = |\mathcal{S}_{xy}|^2$, and hence $|\mathcal{S}_{xy}| = |\mathcal{S}_{x}|$. 
Therefore, we have shown that $|\mathcal{S}_{xy}| = |G|/|X|$ for any $x, y \in X$.
\smallskip

\smallskip

Now, define a set $\mathcal{P}$  by

\begin{equation*}
    \mathcal{P} = \{(g, y) \in G \times X : y \in H \cap gE\}.
\end{equation*}

We can compute the cardinality of $\mathcal{P}$ in two ways using a double-counting argument. On one hand, we have:

\begin{equation}
\label{size of P}
\begin{split}
|\mathcal{P}|&=\sum_{y\in H} \sum\limits_{\substack{g\in G \\ y\in gE}}1=\sum_{y\in H}|\{g\in G: y\in gE\}|\\
&=\sum_{y\in H}\sum_{x\in E}|\{g\in G: y=gx\}|=\sum_{y\in H}\sum_{x\in E}|\mathcal{S}_{xy}|\\
&=\frac{|G||H||E|}{|X|},
\end{split}
\end{equation}
where we have used the fact that $|\mathcal{S}_{xy}|=|G|/|X|$.
On the other hand,

\begin{equation}
\label{upper bound for the size of P}
\begin{split}
|\mathcal{P}|&=\sum_{g\in G}\sum_{y\in H\cap gE}1=\sum_{g\in G} |H\cap gE| \\
    &\leq \max_{g\in G}|H\cap gE|\cdot |G|.    
\end{split}
\end{equation}    

\smallskip

Combining \eqref{size of P} and \eqref{upper bound for the size of P}, one finds that
\begin{equation}
    \max_{g\in G}|H\cap gE|\geq \frac{|H||E|}{|X|},
\end{equation} as claimed in \eqref{lower bound for the intersection}. 
This completes the proof of Theorem \ref{maingroupactiontheorem}. 
\bigskip

\noindent{\bf Proof of  Theorem \ref{BR23}.}
Consider $X = \mathbb{F}_q^d$ and $E \subset \mathbb{F}_q^d$.
For  $r\in (\mathbb{F}_q)^2\setminus \{0\}$, let $\sqrt{r}$ be any square root of $r$, and set
$H = \sqrt{r}E$. Note that $|H|=|E|$.
Finally, let $G$ be the group of translations of $\mathbb{F}_q^d$, acting transitively on $\mathbb{F}_q^d$. 
Theorem \ref{maingroupactiontheorem} implies that
\begin{equation} \label{distancereproveequation}
\max_{a\in \mathbb{F}_q^d}\left|\left(\sqrt{r}E\right) \cap (E + a)\right| \geq \frac{|E|^2}{q^d}.
\end{equation}

If $\lvert E \rvert \geq \sqrt{k+1}q^{\frac{d}{2}}$, then $\max_{a\in \mathbb{F}_q^d} \lvert \sqrt{r}E \cap (E + a) \rvert \geq k + 1$. Thus, there exists an element $a \in \mathbb{F}_q^d$ such that $\lvert \sqrt{r}E \cap (E + a) \rvert \geq k+1$. Consequently, we can establish the existence of a sequence $\{z_1,\dots,z_{k+1}\}$ such that $\{z_1,\dots,z_{k+1}\}\subset \sqrt{r}E \cap (E + a)$. This implies the existence of sequences $\{x_1,\dots,x_{k+1}\} \subset E$ and $\{y_1,\dots,y_{k+1}\} \subset E$, such that $z_i = \sqrt{r}x_i$ and $z_i = y_i + a$ for $1\leq i\leq k+1$.

\smallskip

In summary, we have demonstrated the existence of $(k+1)$-tuples $(x_1, \ldots, x_{k+1}) \in E^{k+1}$ and $(y_1, \ldots, y_{k+1}) \in E^{k+1}$ satisfying the following conditions:

\begin{enumerate}
    \item $x_i \neq x_j$ and $y_i\neq y_j$ for $1 \leq i < j \leq k+1$.

    \smallskip

    \item $y_i + a = \sqrt{r}x_i$ for $i \in \{1, \ldots, k+1\}$. 
\end{enumerate}

\smallskip

Therefore, for $1\leq i<j\leq k+1$, we have: 
\begin{align*}
\begin{split}
    \lVert y_i-y_j\rVert&=\lVert (\sqrt{r}x_i-a)-(\sqrt{r}x_j-a)\rVert=\lVert \sqrt{r}x_i-\sqrt{r}x_j\rVert \\
    &=(\sqrt{r})^2 \lVert x_i-x_j\rVert=r\lVert x_i-x_j\rVert. 
\end{split}
\end{align*}

\smallskip

Since $A$ is a nonempty subset of $\{(i, j) : 1 \leq i < j \leq k+1\}$, we have shown that there exist two $(k+1)$-point configurations, $(x_1, \ldots, x_{k+1}) \in E^{k+1}$ and $(y_1, \ldots, y_{k+1}) \in E^{k+1}$, satisfying the following conditions:

\begin{enumerate} 
\item $x_i \neq x_j$ and $y_i \neq y_j$ for $1 \leq i < j \leq k+1$.

\smallskip

\item $\lVert y_i - y_j \rVert = r \lVert x_i - x_j \rVert$ for $(i, j) \in A$.
\end{enumerate}
This completes the proof of Theorem \ref{BR23}.
\bigskip

\noindent{\bf Proof of  Corollary \ref{sl2actioncorollary}.} Let $G = \textup{SL}_d(\mathbb{F}_q)$, $X={\mathbb F}_q^d$, 
and suppose that $E \subset {\mathbb F}_q^d$ with $|E| \ge \sqrt{k+1}q^{\frac{d}{2}}$. 
The action of $\textup{SL}_d({\mathbb F}_q)$ on ${\mathbb F}_q^d$ is transitive for $d\ge 2$
(this is  just a special case of the general fact, valid for any field $\mathbb F$).
Consequently, Theorem \ref{maingroupactiontheorem} is applicable, and we now show that it leads to the conclusion stated in Corollary \ref{sl2actioncorollary}.
\medskip

If we let $r^{\frac1d}$ be any $d^{\hbox{th}}$ root of $r\in \left(\mathbb F_q\right)^d\setminus \{0\}$,
define $H=r^{\frac{1}{d}}E$.
Then Theorem \ref{maingroupactiontheorem} yields 
$\max_{g\in G}|r^{\frac{1}{d}}E\cap gE|\geq k+1$. 
Let $g\in \textup{SL}_d(\mathbb{F}_q)$ 
be such that $|r^{\frac{1}{d}}E\cap gE|\geq k+1$,
so that there is a set 
$$\{z_1,\dots,z_{k+1}\}\subset r^{\frac{1}{d}}E \cap gE.$$ 

This implies the existence of  $\{x_1,\dots,x_{k+1}\} \subset E$ 
and $\{y_1,\dots,y_{k+1}\} \subset E$, such that $z_i = r^{\frac{1}{d}}y_i$ and $z_i = g x_i$ for $1\leq i\leq k+1$.
\smallskip
 
Hence,
\begin{equation}\nonumber
    \det(z_{i_1},\dots,z_{i_d})=\det(r^{\frac{1}{d}}y_{i_1},\dots,r^{\frac{1}{d}}y_{i_d})=r\det(y_{i_1},\dots,y_{i_d}).
\end{equation}

On the other hand, since $g\in \textup{SL}_d$,
\begin{equation*}
\det(z_{i_1},\dots,z_{i_d})=\det(gx_{i_1},\dots,gx_{i_d})=\det(x_{i_1},\dots,x_{i_d}).
\end{equation*}

Comparing the previous two equalities yields the statement of Corollary \ref{sl2actioncorollary}. 
\medskip

\textbf{Remark.} Let $X=S^{d-1}$ be a $d$-dimensional sphere in $\mathbb{F}_q^d$. 
Then $G=O_d(\mathbb{F}_q)$ acts transitively on $S^{d-1}$. By Theorem \ref{maingroupactiontheorem}
we have that if $E\neq H$ and $|E||H|\geq (k+1)q^{d-1}$, then 
$\max_{g\in O_d}|H\cap gE|\geq k+1$.


\begin{thebibliography}{}


\bibitem{BR23} P. Bhowmik and F. Rakhmonov, \textit{Near optimal thresholds for existence of dilated configurations in $\mathbb{F}_q^d$}, Bulletin of the Australian Mathematical Society. pp. 1-12 (2024)


\bibitem{DORZ23} X. Du, Y. Ou, K. Ren, and R. Zhang, {\it New improvement to Falconer distance set problem in higher dimensions}, (arXiv:2309.04501) (2023). 

\bibitem{DZ19} X. Du and R. Zhang, {\it  Sharp $L^2$ estimate of Schr\"odinger maximal function in higher dimensions}, Ann. of Math. (2) 189 (2019), no. 3, 837-861.

\bibitem{Fal86} K. J. Falconer, {\it On the Hausdorff dimensions of distance sets}, Mathematika \textbf{32} (1986), 206-212.

\bibitem{GIM21} A. Greenleaf, A. Iosevich and S. Mkrtchyan, {\it Existence of similar point configurations in thin subsets of ${\mathbb R}^d$}, Math. Z. \textbf{297} (2021), no. 1-2, 855-865.

\bibitem{GIOW20} L. Guth, A. Iosevich, Y. Ou and H. Wang, {\it On Falconer's distance problem in the plane}, 
Invent. Math. 219 (2020), no. 3, 779-830.

\bibitem{M95} P. Mattila, {\it Geometry of sets and measures in Euclidean spaces}, Cambridge Studies in Adv. Math. {\bf 44}.  Cambridge Univ. Press, 1995.


\bibitem{MR4609035}F. Rakhmonov, {\it Distribution of similar configurations in subsets of $\mathbb{F}_q^d$}, { Discrete Math.} \textbf{346}(10) (2023), Paper No. 113571, 21 pages.


\bibitem{W99} T. Wolff, {\it Decay of circular means of Fourier transforms of measures}, International Mathematics Research Notices \textbf{10} (1999) 547-567.


\end{thebibliography}
\end{document}